\documentclass[a4paper,12pt]{amsart}

\usepackage{amssymb,amscd,amsthm,mathrsfs}
\usepackage[ansinew]{inputenc}
\usepackage[centertags]{amsmath}
\usepackage{epsfig,graphicx}

 \def\RR{{\mathbb R}}  \def\TT{{\mathbb T}}
   
 \def\ZZ{{\mathbb Z}}

\def\cA{\mathcal{A}}

    \def\cW{\mathcal{W}}
\def\cF{\mathcal{F}}   \def\cR{\mathcal{R}}

\newtheorem*{teo*}{Theorem}
\newtheorem*{prop*}{Proposition}
\newtheorem*{cor*}{Corollary}
\newtheorem*{goal*}{Goal}
\newtheorem*{mainteo}{Main Theorem}
\newtheorem*{adendum}{Addendum}

\newtheorem{teo}{Theorem}[section]

\newtheorem{quest}{Question}
\newtheorem{cor}[teo]{Corollary}

\newtheorem{lema}[teo]{Lemma}
\newtheorem{prop}[teo]{Proposition}

\newcommand{\bi}{\begin{itemize}}
\newcommand{\ei}{\end{itemize}}

\theoremstyle{definition}

\theoremstyle{remark}
\newtheorem{obs}[teo]{Remark}

\newcommand{\demo}[1]{\vspace{.05in}{\sc\noindent Proof #1.}}
\newcommand{\dem}{\vspace{.05in}{\sc\noindent Proof.  }}
\newcommand{\lqqd}{\par\hfill {$\Box$} \vspace*{.05in}}

\newcommand{\eps}{\varepsilon}

\newcommand{\en}{\subset}

\DeclareMathOperator{\spectro}{spr}

\author[R. Potrie]{Rafael Potrie}
\address{CMAT, Facultad de Ciencias, Universidad de la Rep\'ublica, Uruguay}
\urladdr{www.cmat.edu.uy/$\sim$rpotrie}
\email{rpotrie@cmat.edu.uy}

\title[Partially hyperbolic diffeomorphisms isotopic to Anosov]{A few remarks on partially hyperbolic diffeomorphisms of $\TT^3$ isotopic to Anosov}
\thanks{The autor was partially supported by CSIC group 618, FCE-3-2011-1-6749 and the Palis-Balzan research project}

\begin{document}
\begin{abstract}
We show that a strong partially hyperbolic diffeomorphism of $\TT^3$ isotopic to Anosov has a unique quasi-attractor. Moreover, we study the entropy of the diffeomorphism restricted to this quasi-attractor. 

\smallskip


\end{abstract}

\maketitle
\section{Introduction}

\subsection{Geometric structures vs. Robust dynamical properties} 

The motivation of this note is a well known example by Ma\~ne (\cite{ManheContributions}) of a deformation of a linear Anosov diffeomorphism of $\TT^3$ which is robustly transitive while it is not hyperbolic. This example is partially hyperbolic with one-dimensional center direction. The idea of the construction is to develop a DA-type perturbation in a small neighborhood of a periodic point to change its index while keeping enough hyperbolicity to guarantee that the iterates of balls eventually behave as if the map were hyperbolic. This is why the fact that the modification is made in a small ball around a periodic points is essential in his proof of robust transitivity (see \cite{BDV,PujSamSurvey}). These examples have been generalized (see for example \cite{BV}) but the fact that the perturbation is local has always been an essential feature. 

From the work of Ma\~ne (\cite{ManheContributions, ManheECL}), and further generalizations to higher dimensions by Bonatti-Diaz-Pujals-Ures (\cite{DPU,BDP}) we know that $C^1$-robust transitivity implies the existence of certain $Df$-invariant geometric structures, namely \emph{dominated splittings} with some uniform behavior of the jacobian in the extremal subbundles\footnote{In dimension 2, the fact that robust transitivity implies Anosov was obtained in \cite{ManheContributions} but the results in \cite{ManheECL} which provide a different proof are important in the results of \cite{DPU,BDP}.}. The results of Bonatti-Crovisier (\cite{BC}) imply that it is not robust transitivity but rather robust chain recurrence which implies such a structure. An advantage of chain-recurrence is that not-being chain-recurrent is a $C^0$-open property which can be quite easily detected while transitivity is less clear. See \cite{BDV,Chab} for surveys on these topics. 

Although in dimension 2 the converse result holds\footnote{The direct implication in dimension 2 gives that a $C^1$-robustly chain-recurrent diffemorphism is Anosov. Again in dimension two, being Anosov is enough to show robust transitivity.}, already in dimension 3 it is easy to see that one can not expect such a result. It is enough to consider the product of a linear Anosov in $\TT^2$ with a (weak) Morse-Smale diffeomorphism of $S^1$ to obtain a strongly partially hyperbolic diffeomorphism of $\TT^3$ which is (robustly) non-chain recurrent since it has hyperbolic attractors and repellers. More sophisticated examples (and much more surprising) have been constructed recently suggesting that this might be more common than what one would expect (see \cite{BG,Shi}). 

Some semilocal mechanisms have been proposed to be added to partial hyperbolicity in order to obtain robust transitivity, to mention a few, one has the well known \emph{blenders} of Bonatti-Diaz (\cite[Chapter 6]{BDV}) or the SH-property proposed by Pujals-Sambarino (see \cite[Chapter 5]{PujSamSurvey}). However, it is not clear if there can be some global topological property which implies robust transitivity under the assumption of partial hyperbolicity or dominated splitting alone. We can formulate this as follows: 

\begin{quest} Is there a manifold $M$ and an isotopy class of diffeomorphisms of $M$ such that the converse to Bonatti-Diaz-Pujals-Ures result holds? 
\end{quest}

Of course this question is interesting only if the isotopy class does admit robustly transitive diffeomorphisms or at least diffeomorphisms satisfying the conclusions of \cite{DPU,BDP}. On the other hand, let us mention that to the present, this author is not aware of an isotopy class of diffeomorphisms in dimension 3 which is known not to admit robustly transitive diffeomorphisms. In \cite{RH} a related question was posed: Is there a partially hyperbolic diffeomorphism with isometric action on the center bundle which is robustly transitive?. For now, the evidence seems to be negative (see \cite{BG,Shi}), but the problem is still open and definitely interesting. 

From the construction of these examples it could be that the following question survives their techniques of construction of examples:

\begin{quest} Is there an isotopy class of diffeomorphisms of a 3-manifold such that strong partial hyperbolicity implies chain-recurrence\footnote{Let us mention that among robustly chain-recurrent partially hyperbolic diffeomorphisms with one dimensional center those which are robustly transitive form a $C^1$-open and dense subset \cite{AbdenurCrovisier}.}?
\end{quest}

Again, the question must be posed to manifolds admitting strong partially hyperbolic diffeomorphisms (see \cite{BoW,BI,HP2}). The natural candidate to have an affirmative answer is the isotopy class of Anosov in $\TT^3$ more in view of the aforementioned examples of Bonatti-Guelman and Y. Shi. In this note we present some mild results in the direction of an affirmative answer, however, the purpose of the note is not only to present these results but mainly to present some arguments which we believe might have some value in understanding this problem and its difficulties. 

Let us mention that obtaining dynamical consequences from $Df$-invariant geometric structures is quite hard in general. However, analog questions in the conservative setting are better understood (\cite{HHUPHdim3,HamUres}). In the conservative setting, the nature of the problems is different though (transitivity is not an issue, it is ergodicity which captures the attention).  

\subsection{Precise statement of the results} 

Given a homeomorphism $g: X \to X$ of a compact metric space $X$ one denotes as $x \dashv y$ if for every $\eps>0$ there exists a $\eps$-pseudo-orbit from $x$ to $y$ (i.e. there exist points $x=z_0, \ldots z_k=y$ with $k\geq 1$ such that $d(z_{i+1},g(z_i))< \eps$). Inside  the set of \emph{chain-recurrent points} (i.e. points such that $x \dashv x$) it is possible to define the following equivalence relation: 

$$ x \sim y \Leftrightarrow x \dashv y \ \ , \ \ y \dashv x $$

Equivalence classes of this relation are called \emph{chain-recurrence classes} and are now quite well understood for ``typical'' diffeomorphisms of a manifold (\cite{BC,Chab}). It is not hard to see that a homeomorphism is \emph{chain-recurrent} (that is, every point of $X$ verifies $x \dashv x$) if and only if there is a unique chain-recurrence class. Moreover, by a well known result of Conley, this is equivalent with the non-existence of a non-trivial open set $U$ such that $g(\overline{U}) \en U$. 

A particularly important type of chain-recurrence class are \emph{quasi-attractors} which by definition are chain-recurrence classes which verify that they posses a basis of neighborhoods $U_n$ such that $g(\overline{U_n})\en U_n$. A \emph{quasi-repeller} is a quasi-attractor for $g^{-1}$. Quasi-attractors always exist, and from the previous discussion one knows that to prove chain-recurrence it is enough to show that given a quasi-attractor $Q$ and a quasi-repeller $R$ one has that $Q \cap R \neq \emptyset$. 

In this paper we restrict to the case of $C^1$-diffeomorphisms of $\TT^3$. Recall that for a diffeomorphism $f: \TT^3 \to \TT^3$ one can define a matrix $A_f \in GL(3,\ZZ)$ given by the way the map $\tilde f (\cdot) - \tilde f(0)$ acts on $\ZZ^3$ where $\tilde f$ denotes any lift of $f$ to $\RR^3$. A diffeomorphism $f: \TT^3 \to \TT^3$ will be said to be \emph{isotopic to Anosov} if the matrix $A_f$ has no eigenvalue of absolute value equal to one. 

As motivated in the introduction the interest will surround the study of (\emph{strongly}) \emph{partially hyperbolic} diffeomorphisms. A diffeomorphism $f: \TT^3 \to \TT^3$ will be said to be \emph{partially hyperbolic}\footnote{Since this is the only notion considered here we remove the word strong from the definition. Beware in comparing with other literature.} if there exists a $Df$-invariant continuous splitting $T\TT^3 = E^s \oplus E^c \oplus E^u$ into one-dimensional subbundles verifying that for some $N>0$ one has that: 

$$ \|Df^N|_{E^s(x)}\| < \|Df^N|_{E^c(x)} \| < \|Df^N|_{E^u(x)}\| \text{ and } \|Df^N|_{E^s(x)}\| < 1 < \|Df^N|_{E^u(x)}\| $$

\begin{mainteo} Let $f: \TT^3 \to \TT^3$ be a partially hyperbolic diffeomorphism isotopic to Anosov, then, $f$ has a unique quasi-attractor and a unique quasi-repeller.
\end{mainteo} 

A similar result was proved in \cite{HP1} for partially hyperbolic diffeomorphisms of 3-dimensional nilmanifolds. In that case the result is optimal in view of the remarkable examples of Y.Shi (\cite{Shi}) mentioned in the introduction, in the case of the isotopy class of Anosov one can still hope that a stronger result holds: Namely, that such diffeomorphisms are chain-recurrent. The idea of the proof is very similar to the one in \cite{HP1} but the use of pseudo-rotations of $\TT^2$ is quite more involved. In fact, the proof of this theorem was the main motivation for \cite{Pot2}. We remark that the idea of relating the unstable holonomy of a partially hyperbolic diffeomorphism with irrational pseudo-rotations of the torus was first used in \cite{McSwiggen} (see also \cite{PaS}) in the other direction, new examples of irrational pseudo-rotations were constructed by the construction of certain partially hyperbolic diffeomorphisms. 

In addition to the main result, some quantitative results are obtained showing some kind of rigidity which would be present in case one constructs an example of non-chain recurrent partially hyperbolic diffeomorphism of $\TT^3$ isotopic to Anosov. This does not exclude the possibility of a hypothetical counterexample\footnote{Indeed, it is reasonable to expect that it is possible to construct a diffeomorphism isotopic to a linear Anosov which has a partially hyperbolic attractor and a partially hyperbolic repeller which are disjoint and satisfy the properties predicted by the Addendum in the global partially hyperbolic setting. The main difficulty lies in understanding how to construct such and example while remaining partially hyperbolic in the \emph{wandering} region.}, but imposes some restrictions. All along this paper we will assume that the partially hyperbolic diffeomorphism $f$ which is isotopic to a linear Anosov $A_f$ verifies that $A_f$ has two eigenvalues of modulus smaller than one. We shall denote as $\lambda_1, \lambda_2, \lambda_3$ the eigenvalues satisfying $|\lambda_1|\leq |\lambda_2|< 1 < |\lambda_3|$. The results admit symmetric statements in the case where $A_f$ has two eigenvalues of modulus larger than $1$ by applying the same results to $f^{-1}$. 

\begin{adendum} Assume that the quasi-attractor $Q$ and the quasi-repeller $R$ of $f$ do not coincide. Then, the entropy of $f|_{Q}$ coincides with the entropy of $A_f$ (i.e. $\log |\lambda_3|=-( \log |\lambda_1| + \log |\lambda_2|)$) while the entropy of $f|_{R}$ equals $-\log |\lambda_1|$.
\end{adendum}

We remark that the result in the addendum follows also from stronger results of J. Yang\footnote{Presented at the 2nd Palis-Balzan Symposium on dynamical systems. Slides are available in the web page of the conference.} which improves related results by Ures (\cite{Ures}) giving dynamical consequences of partial hyperbolicity in the isotopy class of Anosov.    

\subsection{Organization of the paper}
In section \ref{S.Properties} we introduce the notation for the paper as well as introduce some known properties about partially hyperbolic diffeomorphisms isotopic to Anosov in $\TT^3$ and some results on entropy of general diffeomorphisms and maps. Section \ref{S.AttractingRegions} studies some properties on attracting regions and section \ref{S.Uniqueness} completes the proof of the main theorem, in the remark after the proof of the theorem some extensions are indicated without proofs, since the purpose of the note is more to show the idea than to obtain the result we considered better to present a somewhat weaker result for which the proof is more transparent. The addendum is proved in section \ref{S.Entropy}.  

\smallskip
{\bf Acknowledgements: }\textit{I benefited from discussions with C. Bonatti, S.Crovisier, N. Gourmelon and M.Sambarino.}

\section{Some properties}\label{S.Properties}

\subsection{Standing notation} 
In this paper $f: \TT^3 \to \TT^3$ will denote a partially hyperbolic diffeomorphism isotopic to Anosov. 

We consider $\pi: \RR^3 \to \TT^3$ the universal covering map. The map $\tilde f: \RR^3 \to \RR^3$ will denote a lift to the universal cover and $A_f$ will denote its linear part which will be seen both as a matrix which acts in $\RR^3$ and as a diffeomorphism of $\TT^3$. In the whole paper $\lambda_1,\lambda_2, \lambda_3$ will denote the eigenvalues of $A_f$ and it will be assumed that they satisfy $|\lambda_1| \leq |\lambda_2| < 1 < |\lambda_3|$. 

It follows from \cite{PotPHFol} (see also \cite{HP1}) that $|\lambda_1|<|\lambda_2|$ so it is possible to denote as $E^{ss}_A$ and $E^{ws}_A$ the eigenspaces associated to $\lambda_1$ and $\lambda_2$ respectively. The eigenspace asociated to $\lambda_3$ will be denoted as $E^u_A$ and the notation $E^s_A := E^{ss}_A \oplus E^{ws}_A$ will be used. 

\subsection{Invariant foliations} 

It is well known (\cite{HPS}) that there exist $f$-invariant foliations $\cW^s$ and $\cW^u$ tangent to $E^s$ and $E^u$ respectively. These are called respectively strong stable and strong unstable foliations. 

Recently, in this context it was established that dynamical coherence also holds (\cite{BBI2,PotPHFol}) this means that there exist $f$-invariant foliations $\cW^{cs}$ and $\cW^{cu}$ tangent to $E^s\oplus E^c$ and $E^c \oplus E^u$ respectively which by intersection also give rise to a center foliation $\cW^c$ tangent to $E^c$. 

In the universal cover, the lifts of these foliations will be denoted as $\widetilde \cW^\sigma$ ($\sigma=s,c,u,cs,cu$). Given a foliation $\cF$ the leaf through $x$ will be denoted as $\cF(x)$. 

As a part of the proof in \cite{PotPHFol} some important properties were obtained: There is a \emph{global product structure} between $\widetilde \cW^{cs}$ and $\widetilde \cW^u$, this is, for every $x,y\in \RR^3$ we have that $\widetilde \cW^{cs}(x) \cap \widetilde \cW^u(y) \neq \emptyset$ and consists of a unique point. 

Let us write what we have explained in a proposition for further reference:

\begin{prop}[\cite{PotPHFol}]\label{Prop-CoherenceGPS} 
Let $f: \TT^3 \to \TT^3$ be a partially hyperbolic diffeomorphism isotopic to Anosov, then $f$ is dynamically coherent and there exist global product structure between the foliations $\widetilde \cW^{cs}$ and $\widetilde \cW^{u}$ and between $\widetilde \cW^{cu}$ and $\widetilde \cW^{s}$. 
\end{prop}

This also implies that the foliations are \emph{quasi-isometric} (which we shall not define) and it was proved that it further implies that $f$ is \emph{leaf conjugate} to $A_f$ (see \cite{Ham,HP1}). The following version of leaf conjugacy (see \cite{Ham,HP1}) will not be strictly used, but is included here since it might shed light on some technical issues of the proof of the main theorem (see Remark \ref{Remark-Tecnica}). 

\begin{prop}[Corollary 1.5 of \cite{HP1}]\label{Prop-LeafConj} There exists a homeomorphism $\ell: \TT^3 \to \TT^3$ which lifts to a homeomorphism $L: \RR^3 \to \RR^3$ such that $L(\widetilde \cW^{c}(\tilde f(x)))= A_f(E^{ws}_A + L(x))$. In particular, one has that $L(\widetilde \cW^{cs}(x))= E^s_A + L(x)$ and $L(\widetilde \cW^{cu}(x)) = (E^{ws}_A \oplus E^u_A) + L(x)$. 
\end{prop}

\subsection{Semiconjugacy}

It is a classical result that in the isotopy class of Anosov there exists a semiconjugacy to the linear part (see for example \cite[Section 2.3]{PotTesis} for a proof). 

\begin{prop}\label{Prop-Semiconj} There exists $H: \RR^3 \to \RR^3$ continuous and surjective such that $H \circ \tilde f = A_f \circ H$. Moreover, one has that $H(x+\gamma)= H(x) + \gamma$ for every $\gamma \in \ZZ^3$ so, there exists $h: \TT^3 \to \TT^3$ homotopic to the identity such that $h\circ f = A_f \circ h$. 
\end{prop} 

In particular, one has that there exists $K_0>0$ such that $d(H(x),x) \leq K_0$ for every $x \in \RR^3$. From the semiconjugacy property and the fact that the manifolds $\widetilde \cW^\sigma(x)$ ($\sigma=s,c,u$) are properly embedded and quasi-isometric it follows that the following conditions are satisfied (see \cite[Appendix A]{PotPHFol}): 

\begin{itemize}
\item $H(\widetilde \cW^u(x)) = E^u_A + H(x)$ and $H|_{\widetilde \cW^u(x)}: \widetilde \cW^u(x) \to E^u_A + H(x)$ is a homeomorphism.
\item $H(\widetilde \cW^s(x)) \en E^s_A + H(x)$ and $H|_{\widetilde \cW^s(x)}$ is injective. 
\item $H(\widetilde \cW^c(x)) = E^{ws}_A + H(x)$ and $H|_{\widetilde \cW^c(x)} : \widetilde \cW^c(x) \to E^{ws}_A +H(x)$ is a continuous surjective map which may collapse intervals to points. 
\item $H(\widetilde \cW^{cs}(x))= E^s_A + H(x)$. Moreover, the images of two strong stable leaves may intersect but not topologically cross. 
\end{itemize}

See also \cite{Ham,PotPHFol,HP1}. 

\begin{obs}\label{Remark-PointFibers}
Notice that the fact that for every point $x\in \RR^3$ one has that $H^{-1}(H(x))$ is an interval of $\widetilde \cW^c(x)$ implies that the set of points whose preimage is a point has full Lebesgue measure in $\RR^3$. See also \cite{Ures}.  
\end{obs}

\subsection{Measures of maximal entropy}

Ma\~n\'e's example has a unique measure of maximal entropy which is measure theoretically equivalent to Lebesgue for $A_f$ (\cite{BFSV}). This was generalized by Ures (\cite{Ures,HP1}) for any partially hyperbolic diffeomorphism isotopic to Anosov in $\TT^3$. 

We refer the reader to \cite[Chapters 3 and 4]{KH} or \cite{ManheLibro} for definitions and basic result about measure theoretical and topological entropy. The variational principle (\cite[Chapter IV]{ManheLibro}) implies that the topological entropy equals the supremum of the measure theoretical entropies of all ergodic invariant measures. When a measure $\mu$ verifies that $h_\mu(f)= h_{top}(f)$ it is called a \emph{measure of maximal entropy}. This notion makes sense for compact $f$-invariant subsets which also have a well defined entropy and so the concept makes sense. 

The following holds: 

\begin{teo}\label{Teo-Entropia} For every compact $f$-invariant set $\Lambda$ there exists at least one measure $\mu$ supported on $\Lambda$ such that $h_\mu(f) = h_{top}(f|_\Lambda)$. Moreover, there exists a unique measure $\mu_{max}$ such that $h_{\mu_{max}}(f)=h_{top}(f)= h_{top}(A_f)=\log |\lambda_3|$.  
\end{teo}

The first statement follows directly of applying \cite[Corollary 1.3]{DFPV} to $\varphi=1$ (see also \cite{CoYo}). The rest follows from \cite{Ures,HP1}. 

\subsection{The entropy conjecture in the first homology group}

We state here a well known result of Manning \cite{Manning}, but first some definitions must be introduced. 

Let $T: M \to M$ be a continuous map of a compact manifold $M$ (possibly with boundary). Denote as $H_1(M)$ the first homology group with real coefficients which is finite dimensional since $M$ is compact. The map $T$ induces a map $T_\ast : H_1(M) \to H_1(M)$ which can be represented as a real matrix. We define $\spectro(T_\ast)$ to be the modulus of the eigenvalue of $T_\ast$ of largest modulus.  

Denote as $\Lambda = \bigcap_{n\geq 1} T^n(M)$.

\begin{teo}[Corollary 8.1.3 of \cite{KH}]\label{TEo-ConjEntro} For $T: M \to M$ as above it holds $h_{top}(T) \geq \log \spectro (T_\ast)$.
\end{teo}

\section{Attracting regions}\label{S.AttractingRegions} 

If $f$ is chain-recurrent, the Main Theorem holds trivially being the whole $\TT^3$ both the unique quasi-attractor and quasi-repeller. 

Along this section we shall study consequences of the existence of an open set $U$ such that $\emptyset \neq U \neq \TT^3$ verifying that $f(\overline{U})\en U$. Recall that the non-existence of such $U$ is equivalent to chain-recurrence of $f$. (This is a consequence of Conley's theorem, see for example \cite{Chab}.)

Since $d(U^c, f(\overline{U})) > a>0$, without loss of generality, one can assume that $\overline{U}$ is a compact manifold with boundary as well as $U^c$. 

Let us denote as $\cA = \bigcap_{n\geq 0} f^n(\overline U)$ and $\cR = \bigcap_{n\geq 0} f^{-n}(U^c)$. It is clear that $\cA$ and $\cR$ are non-empty disjoint compact $f$-invariant sets.  

\begin{lema}\label{Lema-SaturadosPorFuertes} For every $x\in \cA$ one has that $\overline{\cW^u(x)} \en \cA$. Symmetrically, for every $x\in \cR$ one has that $\overline{\cW^s(x)}\en \cR$. 
\end{lema}    

Notice that this lemma, whose proof is left to the reader, also applies to quasi-attractors and quasi-repellers which are (at most countable) intersections of such kind of sets. 

Denote as $\imath_1 : \overline{U} \to \TT^3$ and $\imath_2: U^c \to \TT^3$ the inclusion maps. These maps induce group morphisms $(\imath_1)_\ast: \pi_1(\overline{U}) \to \pi_1(\TT^3)$ and $(\imath_2)_\ast: \pi_1(U^c) \to \pi_1(\TT^3)$. Define $\Gamma_1= (\imath_1)_\ast(\pi_1(\overline{U}))$ and $\Gamma_2 = (\imath_2)_\ast(\pi_1(U^c))$. One can identify $\pi_1(\TT^3)$ with $\ZZ^3$ and so $A_f$ acts in $\ZZ^3$ as the action of $f$ in $\pi_1(\TT^3)$.

\begin{lema}\label{Lema-Gama1YGama2NoTrivialesEInvariantes} The groups $\Gamma_1$ and $\Gamma_2$ are $A_f$-invariant and non-trivial.
\end{lema}

\dem The proof will be given for $\Gamma_1$. An analogous argument applied to $f^{-1}$ gives the respective result for $\Gamma_2$. 

First, we show that $\Gamma_1$ is not trivial. To do this, consider the lift $\tilde U$ of $U$ and $U_0$ a connected component of $\tilde U$ such that there is a point $x\in U_0$ such that $\pi(x) \in \cA$. Since $\cA$ is $\cW^u$-saturated and $d(\cA, \partial U)>a>0$ one has that the $a$-neighborhood of $\widetilde \cW^u(x)$ is contained in $U_0$. Since $\widetilde \cW^u$ is properly embedded, it follows that the volume of $U_0$ is infinite, in particular, there must exist a deck transformation $\gamma \in \ZZ^3\setminus \{0\}$ such that $U_0 \cap (U_0 +\gamma) \neq \emptyset$, but this implies that $U_0 = U_0 + \gamma$. This in turn implies that $\gamma \in \Gamma_1$ showing that $\Gamma_1$ is not trivial. 

To show invariance, consider $\eta \in U$ a loop representing an element $\gamma \in \Gamma_1$. It follows that $f(\eta)$ represents $A_f(\gamma)$ and since $f(\overline U)\en U$ it follows that $A_f(\gamma) \in \Gamma_1$. This shows that $A_f(\Gamma_1) \en \Gamma_1$. 

We claim that this implies that $A_f(\Gamma_1)= \Gamma_1$. Indeed, consider the increasing\footnote{Since $A_f(\Gamma_1)\en \Gamma_1$ it follows that if $m < k$ then $A_f^{-m}(\Gamma_1) \en A_f^k(\Gamma_1)$.} union $\tilde \Gamma_1= \bigcup_{n\geq 0} A_f^{-n}(\Gamma_1)$ and let $S \en \tilde \Gamma_1$ be a finite set which generates $\tilde \Gamma_1$. There exists $n_0$ such that $A_f^{-n_0}(\Gamma_1)$ contains $S$. This implies that $A_f^{-1}(A_f^{-n_0}(\Gamma_1)) = A_f^{-n_0}(\Gamma_1)$ and since $A_f$ is an isomorphism $A_f^{-1}(\Gamma_1)= \Gamma_1$. 

\lqqd

\begin{cor}\label{Coro-Gama1YGamma2SonTodo} Both $\Gamma_1$ and $\Gamma_2$ have rank 3 as subgroups of $\ZZ^3$. 
\end{cor} 

Via the Hurewicz morphism (see for example \cite[Chapter 3.1.e]{KH}) the induced map in homology (with real coeficients) which we denote also as $(\imath_1)_\ast: H_1(\overline{U}) \to H_1(\TT^3)$ is linear and surjective. The same holds for $(\imath_2)_\ast$.  


\begin{prop}\label{Prop-Entropy} The entropy of $f|_{\cA}$ is equal to $h_{top}(f)=\log|\lambda_3|$. The entropy of $f|_{\cR}$ is contained in $[-\log|\lambda_1| , \log|\lambda_3|)$. 
\end{prop}

\dem First we prove the result for $f|_{\cA}$. Notice that we have the following commuting diagram:

\[\begin{array}[c]{ccc}
H_1(\overline U)&\stackrel{(f|_{\overline{U}})_\ast}{\rightarrow}&H_1(\overline{U})\\
\downarrow\scriptstyle{(\imath_1)_\ast}&&\downarrow\scriptstyle{(\imath_1)_\ast}\\
H_1(\TT^3)&\stackrel{f_\ast}{\rightarrow}& H_1(\TT^3)
\end{array}\]

\noindent which clearly implies that $|\spectro((f|_{\overline U})_\ast)|
\geq |\spectro (A_f)|$. 

%
%
%
%
%
%
%
%

Since we have assumed that $\overline U$ is a compact manifold with boundary, by Theorem \ref{TEo-ConjEntro} we obtain that $h_{top}(f|_{\cA}) \geq \log |\spectro (A_f)|= \log |\lambda_3|$. Since $h_{top}(f) \leq \log |\lambda_3|$ by Theorem \ref{Teo-Entropia} we get the desired result. 

The same argument applied to $U^c$ and $f^{-1}$ gives that $h_{top}(f|_{\cR}) \geq \log |\spectro(A_f^{-1})| = -\log |\lambda_1|$. On the other hand, since $\cR$ is compact and disjoint with $\cA$ and admits a measure realizing the entropy $h_{top}(f|_{\cR})$ by Theorem \ref{Teo-Entropia} it follows that the entropy cannot be equal to $\log |\lambda_3|$ by uniqueness of the measure of maximal entropy.
\lqqd

\begin{obs}\label{Obs-hAesTodo} The asymmetry in the statement reflects the fact that $A_f$ has two stable eigenvalues and one unstable eigenvalue. Since $h(\cA)$ is $\cW^u$-saturated and compact one has that $h(\cA)=\TT^3$. Since there must be points such that $h^{-1}(x)$ is a unique point, this cannot hold for $h(\cR)$. 
\end{obs}

\section{Uniqueness of attractors and repellers}\label{S.Uniqueness} 

This section shows that there exist a unique quasi-attractor and a unique quasi-repeller for $f$. This will complete the proof of the main Theorem. An important remark is that the argument is not symmetric since $f$ and $f^{-1}$ are essentially different as it was revealed in Proposition \ref{Prop-Entropy}. In fact, uniqueness of the quasi-attractor is almost direct:

\begin{lema}\label{lema-UniqueAttractor} There exist a unique quasi-attractor for $f$.
\end{lema}

\dem Assume there are two quasi-attractors $Q_1$ and $Q_2$. Since they are both compact and $\cW^u$-saturated (Lemma \ref{Lema-SaturadosPorFuertes}) one gets that, as in Remark \ref{Obs-hAesTodo}, that $h(Q_1)=h(Q_2)= \TT^3$. On the other hand, since there are points in $\TT^3$ such that $h^{-1}(x)$ is a single point (Remark \ref{Remark-PointFibers}) it follows from this that $Q_1 \cap Q_2 \neq \emptyset$. Since they are chain-recurrence classes this implies that $Q_1=Q_2$ giving uniqueness. 
\lqqd

The proof of the uniqueness of the quasi-repeller is by contradiction and will occupy the rest of this section. 


\begin{obs}\label{Remark-Tecnica} The proof would be more direct if one could show the existence of a torus $T$ (topologically) transverse to the unstable foliation $\cW^u$ in $\TT^3$ and show that the return map to that torus is a non-resonant torus homeomorphism which is as it is done in \cite{HP1} (where a result similar to Proposition \ref{Prop-LeafConj} above in that context is used). In \cite{HP1} once one obtains such torus homeomorphism, one applies techniques about such homeomorphisms to obtain the result. In the context of \cite{HP1} (nilmanifolds) though, it is possible to conclude in a direct way. Here the argument is slightly more delicate (in particular, by contradiction) and uses the results of \cite{Pot2}. Moreover, the existence of such tori is less clear, so we argue in a slightly different way to get the non-resonant torus homeomorphism.  
\end{obs}

\demo{of the Main Theorem}
Assume that there exist two disjoint quasi-repellers $R_1$ and $R_2$ and consider open sets $U_1$ and $U_2$ such that $f(\overline{U_i}) \en U_i$ ($i=1,2$) and such that $R_1 \en U_1^c$ and $R_2 \en U_2^c$. If $R_1 \neq R_2$ it is possible to choose $U_1$ and $U_2$ such that $R_2 \en U_1$. Modulo considering an iterate and modifying $U_1$ it is possible to assume that $U_1$ is connected. 

We consider, as in the previous section, the sets $\cA= \bigcap_{n\geq  0} f^n(\overline{U_1})$ and $\cR = \bigcap_{n\geq 0}f^{-n}(U_1^c)$. From the choice of $U_1$ and invariance of $R_1$ and $R_2$ it follows that: 

$$  R_1 \en \cR \text{ and } R_2 \en \cA $$

\noindent we will show that this implies that $\cA$ intersects $\mathcal{R}$, a contradiction.  Notice that by the assumption on $U_1$ one has that $\cA$ is connected. 

Consider a point $x_0 \in \RR^3$, using the global product structure (Proposition \ref{Prop-CoherenceGPS}) we can define a $\ZZ^3$-action on $\widetilde \cW^{cs}(x_0) \cong \RR^2$ as follows: 
$$ \alpha : \ZZ^3 \times \widetilde \cW^{cs}(x_0) \to \widetilde \cW^{cs}(x_0) \ ; \ (\gamma, y) \mapsto \widetilde \cW^{cs}(x_0) \cap \widetilde \cW^{u}(y + \gamma) $$
Using the semiconjugacy $H$ it is clear that this action is semiconjugated with the linear action: 
$$ \beta : \ZZ^3 \times E^s_A \to E^s_A  \ ; \ (\gamma,z) \mapsto E^s_A \cap (E^u_A + z + \gamma) $$
Since $\beta$ induces an irrational translation on $\TT^2$ by making the quotient by a $\ZZ^2$ subgroup of $\ZZ^3$ it follows that $\alpha$ defines a non-resonant homeomorphism\footnote{This means that it has a unique rotation vector which is totally irrational, see \cite{Pot2}.} of $\TT^2$. More precisely, this means that there exists a subgroup $\Gamma \en \ZZ^3$ isomorphic to $\ZZ^2$ such that under the action of $\alpha$,  $\widetilde \cW^{cs}(x_0)/\alpha(\Gamma)$ is homeomorphic to $\TT^2$ and the action of $\alpha(\gamma)$ for some $\gamma \notin \Gamma$ one has that $\alpha(\gamma)$ is the lift of a non-resonant torus homeomorphism. Let us call $F$ to the homeomorphism of $\TT^2 \cong \widetilde \cW^{cs}(x_0)/\alpha(\Gamma)$ which lifts to $\alpha(\gamma)$ under the covering $\tilde \pi : \widetilde \cW^{cs}(x_0) \mapsto \widetilde \cW^{cs}(x_0)/\alpha(\Gamma)$.  

In \cite{Pot2} the following result is proved:

%


\begin{prop}[Proposition B \cite{Pot2}]\label{PropositionB} 
Given a compact connected set $\Lambda$ such that $F(\Lambda) \en \Lambda$ one has that for every neighborhood $V$ of $\Lambda$ there exists $K>0$ such that every connected set $\Lambda'$ having diameter larger than $K$ in the universal cover intersects $V$. 
\end{prop} 

Since $\cA$ is $\cW^u$-saturated and contains an entire $\cW^s$-leaf (because $R_2 \en \cA$) one gets that the intersection of $\tilde \cA$ with $\widetilde \cW^{cs}(x_0)$ contains a non-bounded connected set (the intersection of the saturation by $\widetilde \cW^u$ of $\tilde R_2$ with $\widetilde \cW^{cs}(x_0)$) which is invariant under $\alpha(\gamma)$. This implies that for $F$ there is a closed connected set which is forward invariant. 

Applying the previous proposition one deduces that the quasi-repeller $R_1$ must intersect every neighborhood of $\cA$ and since both are closed sets, they must intersect. This gives a contradiction and concludes the proof of the main theorem. 

\lqqd

\begin{obs}
\begin{itemize}
\item[(a)] The same proof can be used to show that both the quasi-attractor and the quasi-repeller are connected. In fact, it can be proved that there is a unique minimal set of each strong foliation (but the argument strongly depends on the one-dimensionality of the center direction).   
\item[(b)] In \cite{PotPHFol} the case where the splitting is of the form $T\TT^3= E^{cs} \oplus E^u$ was studied under a hypothesis which was called \emph{almost dynamical coherence}. The same strategy can be used in that case but weaker results are obtained.
\begin{itemize}
\item[--] If $E^{cs}$ is volume contracting (see \cite{BDP}) then, the argument implies that there is a unique quasi-repeller\footnote{ The case where $A_f$ has two eigenvalues smaller than one is considerably easier and very similar to the proof above, but the other case is a little bit more complicated since the induced dynamics is no longer semiconjugate to a rigid translation.}.
\item[--] In any case (without any assumptions on the behavior on $E^{cs}$) it is possible to show by similar arguments that for $C^1$-generic diffeomorphisms homoclinic classes have empty interior (see \cite{ABD,PotS}). 
\item[--] The argument of the ``easy'' Lemma \ref{lema-UniqueAttractor} fails short of giving a unique quasi-attractor\footnote{The question of wether a non-resonant torus homeomorphism may have more than one minimal set is, I believe, open. Even in the case where it is semiconjugate to a rigid translation.} in the case $T\TT^3 = E^{cs} \oplus E^u$. 
\end{itemize}
\end{itemize}
\end{obs}

\section{Entropy on the repeller}\label{S.Entropy} 

In this section we prove the Addendum to the main theorem. It remains to show that the entropy restricted to $\cR$ cannot be larger than $-\log |\lambda_1|$, this is proved in the following Lemma.  The proof of this lemma relies on ideas from \cite{BFSV,HSX,Ures} so some familiarity with those papers will be assumed (precise references will be given). 

\begin{lema} The entropy of $f|_{\cR}$ equals $-\log|\lambda_1|$.
\end{lema}

\dem We must show that $h_{top}(f|_{\cR}) \leq -\log |\lambda_1|$ since the other inequality was established in Proposition \ref{Prop-Entropy}. 

To do this, consider a measure $\nu$ supported in $\cR$. Theorem 3.3 in \cite{HSX} applied to $f^{-1}$ implies that $h_\nu(f^{-1}) \leq \lambda^c(\nu) + |\log|\lambda_1||$. This is because the volume expansion along the strong stable foliation is bounded by $-\log|\lambda_1|$ due to quasi-isometry of the strong stable foliation (see the proof of Theorem 5.1 in \cite{Ures}). 

In principle, it may hold that the center-Lyapunov exponent of $\nu$ is larger than zero. 

We will show that in this case there exists a measure $\nu'$ supported in $\cR$ with the same entropy as $\nu$ verifying that $\lambda^c(\nu') \leq 0$ for $f^{-1}$. 

To do this, notice that $h_\nu(f^{-1})$ is the same as $h_{h^\ast(\nu)}(A_f)$ where $h^\ast$ denotes the push-forward of the measure. This follows from the same argument in \cite{BFSV,Ures} using the fact that fibers of the semiconjugacy are bounded intervals and thus they carry no entropy. 

Consider the measure $\nu'$ supported on $\cR$ given by the lift of $h^\ast(\nu)$ to $\TT^3$ such that its desintegration along $\cW^c$ consists of the sum of two dirac measures of $1/2$ of weight in each extremal point of $h^{-1}(x) \cap \cR$. This measure is clearly invariant and supported on $\cR$ and since $\cR$ is a topological repeller it follows that the center exponent is $\leq 0$ for $f^{-1}$. 

\lqqd

%
%
%
%



\end{document}